\documentclass[11pt,leqno]{article}
\usepackage[all]{xy}
\usepackage{amssymb,amsmath, amsthm}
\theoremstyle{definition}

\theoremstyle{definition}

\theoremstyle{remark}

\usepackage{amssymb, palatino}
\usepackage{graphicx}

\begin{document}
\title{On the fundamental groups of compact Sasakian manifolds  }
\author{Xiaoyang Chen}
\date{}
\maketitle
\begin{abstract}
We study the fundamental groups of compact Sasakian manifolds, which we call Sasaki groups. It is shown that all known K$\ddot{a}$hler groups are Sasaki, in particular, all finite groups are Sasaki. On the other hand, we show there exists many restrictions on the fundamental groups of compact Sasakian manifolds. We also study the Abel-Jacobi map of a compact Sasakian manifold and its applications to Sasaki groups.
\end{abstract}
\section*{1. Introduction}
 Recently Sasakian manifolds attract a lot of attention. See [BoG08] for a comprehensive introduction.  However, very little is
known for the fundamental groups of compact Sasakian manifolds.
 Here we make some steps on this direction.
\par Sasakian manifolds are odd dimensional analogue of K$\ddot{a}$hler manifolds. Recall a Riemannian manifold $(M^{2n+1},g)$ ( We always consider connected manifolds in this paper) to be Sasakian if it has a unit Killing field $\xi$,
satisfying the equation
$$ R(X,\xi)Y=g<\xi,Y>X-g<X,Y>\xi$$
Given such a characteristic vector $\xi$ (also called Reeb vector field), we define a $(1,1)$ tensor $\phi$ to be given by $\phi(X)=\nabla_X \xi $ and the characteristic one form $\eta$
to be given by $\eta(X)=<X,\xi>.$ Altogether we call $(g, \xi, \eta, \phi)$ a Sasakian structure.
The vector field $\xi$ defines the characteristic foliation $F_{\xi}$ with one dimensional leaves and the kernel of $\eta$ denoted by $D$, called the contact bundle,
inherits an almost complex structure by restriction of $\phi$.
Let $g^T=g-\eta \otimes \eta$.
It turns out $(g, \xi, \eta, \phi)$ is a Sasakian structure iff $(D, g^T, \phi|_D, d\eta)$ defines a transversal K$\ddot{a}$hler structure with transversal K$\ddot{a}$hler form $d\eta$.
\par A Sasakian structure on $M$ is called quasi regular if all leaves of characteristic foliation $F_{\xi}$ are closed, otherwise it is called irregular. By a theorem of Wadsley [Wa75], if a Sasakian structure on $M$ is quasi regular, $\xi$
generates a locally free $S^1$ action on $M$.
A Sasakian structure on $M$ is called regular if this action is free.
If the Sasakian structure on $M$ is quasi regular, then the quotient space $M/F_{\xi}$ is a K$\ddot{a}$hler orbifold. In gerenal, there is no quotient space if it is irregular.
\par There is a natural transversal Levi-Civita connection $\nabla^T$ on $M$ defined by
$\nabla_X^T Y=[\nabla_X Y]^p$ if $X,Y\in D$
and $\nabla_\xi^T Y=[\xi, Y]^p$, $\nabla_X ^T \xi=0$, where we denote $Z^p$ is the projection of $Z$ to $D$ for any $Z\in TM$ and
$\nabla$ is the usual Levi-Civita connection induced by $g$.
\par Now define the transversal curvature tensor $R^T$ by
$$ R^T(X,Y)Z=\nabla_X^T\nabla_Y^TZ-\nabla_Y^T\nabla_X^TZ-\nabla^T_{[X,Y]}Z.$$
$$ R^T(X,Y,Z,W)=g<R^T(X,Y)Z,W>.$$
and the transversal sectional curvature of the plane spanned by $X,Y$ is defined by
$$ K^T(X,Y)=R^T(X,Y,Y,X)$$
\\where $X,Y, Z, W \in D$.
\par We say $(M,g)$ has nonpositive transversal sectional curvature if and only if $K^T(X,Y)\leq 0$ for any linearly independent vectors $X, Y \in D$. Examples of compact Sasakian manifolds
with nonpositive transversal sectional curvature are given by Boothby-Wang fibrations over compact locally Hermitian symmetry spaces of noncompact type. See section 2 for details on Boothby-Wang fibration.
\par  We say a finitely generated group is Sasaki (K$\ddot{a}$hler) if it is the the fundamental group of a compact Sasakian (K$\ddot{a}$hler) manifold.
\par It is known by various geometers that any compact Sasakian manifold has even first Betti number, see
[BG67], [F66], [Tan67]. Since one can lift a Sasakian structure to the covering space of a Sasakian manifold, the following proposition is immediate.
\\
\\
\noindent \textbf{Propositon 1.1} {\itshape  A group which contains a finite index subgroup with odd first Betti number is not Sasaki. In particular, nontrivial free groups can not be Sasaki. }
\\
\par The next proposition gives many examples of Sasaki groups.
\\
\\
\noindent \textbf{Propositon 1.2} {\itshape  If $\Gamma$ is the fundamental group of some compact Hodge manifold, it is also Sasaki. In partiular, all known K$\ddot{a}$hler groups are Sasaki.}
 \\
 \par Using a classic result of Serre ([ABCKT96]) that any finite group is the fundamental group of some compact Hodge manifold
 and also propositions 1.1, 1.2, one easily gets
 \\
 \\
 \noindent \textbf{Corollary 1.3}  {\itshape All finite groups are Sasaki and an abelian group
 is Sasaki if and only if it has even rank.}
 \\
 \par It was conjectured that any  K$\ddot{a}$hler group is the fundamental group of some compact Hodge manifold.
 So it is natural to propose the following
 \\
 \\
 \noindent \textbf{Conjecture} {\itshape  All K$\ddot{a}$hler groups are Sasaki.}
 \\
\par It is easy to see there exists Sasaki groups which are not K$\ddot{a}$hler. For example, the discrete,
 torsion free and cocompact subgroups of real Heisenberg group $H^{2n+1}$ are Sasaki. However,
they are not K$\ddot{a}$hler if $n\leq 3$ by a theorem of Carlson and Toledo, see [CaT95].
\\
 \par The next theorem tells us some well known results for  K$\ddot{a}$hler groups remain true for Sasaki groups.
  \\
  \\
  \noindent \textbf{Theorem 1.4}
  {\itshape Suppose $\Gamma$ is a Sasaki group, then
  \\ (1) $\Gamma$ has either $0$ or $1$ end; in particular, $\Gamma$ can not split as a nontrivial free product;
  \\ (2) If $\Gamma$ is solvable, it
   contains a nilpotent subgroup of finite index.}
  \\
  \par The following proposition gives more restrictions on the fundamental groups of compact regular Sasakian manifolds.
  \\
  \\
 \noindent \textbf{Proposition 1.5 } {\itshape Suppose $\Gamma$ is the fundamental group of some compact regular Sasakian manifold. Then
 \\ (1) $\Gamma$ can not be a cocompact, discrete and torsion free subgroup of $G$, where $G$ is $SO(1,n), n>2$ or $F_{4(-20)}$
or a simple real non-Hermitian Lie group of noncompact type with real rank at least 20.
\\
 (2) $\Gamma$ is the fundamental group of some compact
 $3$-manifold $V^3$ if and only if the universal cover of $V^3$ is differemorphic to $S^3$ or the $3$-dimensional Heisenberg group or $\widetilde{SL(2,\mathbb{R})}$.}
\\
\par It is quite believable that proposition 1.5 holds for all Sasaki groups without any regular assumption.
\\
\\
 \noindent \textbf{Theorem 1.6 } {\itshape Suppose $(M^{2n+1}, g)$ is a compact Sasakian manifold with nonpositive transversal sectional curvature, then $\pi_1(M^{2n+1})$ is an infinite group which can not be realized as the fundamental group of any compact Riemannian manifold with nonpositive sectional curvature. }
 \\
 \par  As a consequence of theorem 1.6, we recover the following fact which was previous
  known by P. Eberlein's splitting theorem ([Eb82]): If $\Gamma = \pi_1(M^3)$, where
 $M^3$ is a compact 3-manifold with geometry modeled on $\widetilde{SL(2,\mathbb{R})}$, $\Gamma$ can not be
 the fundamental group of any compact Riemannian manifold with nonpositive sectional curvature. However, as far as the author konws, theorem 1.6 does not follow from P. Eberlein's theorem if $n>1$.
 \par Theorem 1.6 suggests to study a new class of infinite groups, which we hope to cover in the future.
 \\
 \\
 \par The organization of this paper is the following: In section 2 we review two constructions in Sasakian geometry and prove proposition 1.2. Theorem 1.4 is proved in section 3. In section 4 we study harmonic maps from compact Sasakian manifolds and prove proposition 1.5. In section 5 we prove theorem 1.6. In section 6 we study the Abel-Jacobi map of a compact Sasakian manifold and its applications to Sasaki groups.
\\
\\
\noindent \textbf{Acknowledgement:} The author would like to express his gratitude to
 Professor Karsten Grove, who is his advisor, for insistent encouragement and helpful advice.
  He is also grateful to
 Professor Jianguo Cao, who passed away during the writing of this paper.
 He also thanks
 Professor Marc Burger, Frederic Campana,  Qintao Deng, Isaac Goldbring,
 Wolfgang Heil, Michael Kapovich, Shmuel Weinbrger, Gaofeng Zheng and others for many helpful discussions.
He is also grateful to Weimin Peng and his family for insistent encouragement.

\section*{2. Constructions in Sasakian geometry and proof of proposition 1.2}
We first recall two basic constructions in Sasakian geometry. For details, see [BoG08]. For simplicity, we restrict ourselves to the class of regular Sasakian manifolds.
\\
\\
\noindent \textbf{ Boothby-Wang fibration:} Suppose $(N^{2n}, \omega)$ is a compact Hodge manifold with
integral K$\ddot{a}$hler form $\omega$. By a theorem of Kobayashi [Kob63], there exists
 a principal circle bundle $P$ over $N^{2n}$ and a connection form $\eta$ on $P$ such that $d\eta =p^*\omega$, where
  $p$ is the projection map. Let $M^{2n+1}$ be the total space of $P$.
Now define a Riemanian metric $g$ on $M^{2n+1}$ by $g=p^*h+ \eta \otimes \eta$, where $h$ is the associated K$\ddot{a}$hler metric for $\omega$ on $N^{2n}$. It is not hard to check that $(M^{2n+1}, g)$ becomes a Sasakian manifold with transversal K$\ddot{a}$hler form $\omega$. We call $(M^{2n+1}, g)$ is the Boothby-Wang fibration over $(N^{2n}, h)$.
It is not hard to extend this construction to the case when $(N^{2n}, \omega)$ is a compact Hodge orbifold.
\par From the construction of Sasakian structure on $M^{2n+1}$, it is easy to see $(M^{2n+1}, g)$ has nonpositive transversal sectional curvature if $(N^{2n}, h)$ has nonpositive sectional curvature, for example, $(N^{2n}, h)$ is a compact locally Hermitian symmetry space of noncompact type.
\\
\\
\noindent \textbf{ Join construction: } Suppose $M_1$ and $M_2$ are two compact Sasakian manifolds
over compact Hodge manifolds $(N_1, \omega_1)$ and $(N_2, \omega_2)$, respectively. Then $(N_1 \times N_2, \omega_1 + \omega_2)$ is also a compact Hodge manifold. Denote $M_1* M_2$ be the Sasakian manifold over $N_1 \times N_2$ coming from Boothby-Wang fibration. We call $M_1* M_2$ is the join of $M_1$ and $M_2$. It turns out $M_1* M_2$ is a $M_2$ bundle over $N_1$, see [BoG08]. It is easy to see compact Sasakian manifolds
with  nonpositive transversal sectional curvature
are closed under join construction.
\\
\\
\noindent \textbf{Proof of proposition 1.2}
\par
Suppose $\Gamma=\pi_1(N)$, where $N$ is a compact Hodge manifold. Let $M$ be the Sasakian manifold over $N$ coming from Boothby-Wang fibration.
 Let $V=M * S^3$. Then $V$ is a compact Sasakian manifold and also a $S^3$ bundle over $N$. From the long exact sequence of homotopy groups, one gets $\pi_1(V)\simeq \Gamma$.
\section*{3. Orbifold fundamental groups of compact K$\ddot{a}$hler orbifolds}
\par We refer the readers to [Da08] or [Th78]  for general notions on orbifolds and orbifold fundamental groups.
\par In [C11], the author studied the orbifold fundamental groups of K$\ddot{a}$hler orbifolds. In particular, he proved the following:
\\
\\
\noindent \textbf{Lemma 3.1} {\itshape Suppose $G$ is the orbifold fundamental group of a compact K$\ddot{a}$hler orbifold, then
\\ (1) $G$ has $0$ or $1$ end;
\\ (2) If $G$ is solvable, it contains a nilpotent subgroup of finite index.}
\\
\\
\noindent \textbf{Remark:} The first part of lemma 3.1 is not explicitly stated in [C11], however, it follows from the arguments there and the corresponding result in the manifold case.
\\
\par We recall some basic facts on the end of a finitely generated group. See [DK] for more details.
Suppose $X$ is a locally compact connected topological space. The set of ends of $X$, denoted by $E(X)$, is defined as the inverse limit: \[\lim_{K\subseteq X} \pi_0(K^c)\], where $K$ is a compact subset of $X$, $K^c$ is the compliment of $K$ and $\pi_0(K^c)$ is the number of path-connected components of $K^c$. If $\Gamma$ is a finitely generated group, the space of ends $E(\Gamma)$ is defined as the set of ends of its Cayley graph. The elements in $E(\Gamma)$ are called the ends of $\Gamma$. We denote $ e(\Gamma)$ be the cardinality of
  $E(\Gamma)$. It can be shown $e(\Gamma)$ is a quasi-isometric invariant of $\Gamma$.
The following facts are due to Hopf and Freudenthal:
\\
\par 1. $ e(\Gamma)\in \{0,1,2,\infty\}.$
\par 2. $e(\Gamma)=0$ if and only if $\Gamma$ is finite.
\par 3. $e(\Gamma)=2 $ if and only if $\Gamma$ has an infinite cyclic subgroup of finite index.
\\
\par It is a well known theorem of Stalling that a finitely generated group has more than one end if and only if
it is a nontrivial amalgamated free product or an HNN extension over a finite subgroup, see [Be68], [St68], [St71].
\\
\\
\noindent \textbf{Lemma 3.2} {\itshape Suppose we have the following short exact squence of groups
$$1 \rightarrow \mathbb{Z} \rightarrow A_1 \rightarrow A_2 \rightarrow 1$$
then $A_1$ has $1$ or $2$ ends.}
\begin{proof}
If $A_2$ is finite, $A_1$ contains an infinite cyclic subgroup of finite index and hence has $2$ ends.
Otherwise $A_2$ is infinite, it follows that $A_1$ has $1$ end by proposition 1.9 in [Co70].
\end{proof}
\noindent \textbf{Lemma 3.3} {\itshape Suppose we have the following short exact squence of groups
$$1 \rightarrow B_1 \rightarrow B_2 \rightarrow B_3 \rightarrow 1$$
and $B_1$ is a cyclic group and $B_3$ is virtually nilpotent, then $B_2$ is virtually nilpotent.
Here we say a group is virtually nilpotent if it contains a nilpotent subgroup of finite index.}
\begin{proof}
By passing to a subgroup of finite index, we can assume $B_3$ is a nilpotent group. Since $B_1$ is normal in $B_2$, we see $B_2$ acts on $B_1$ by conjugation. This gives a group homomorphism from $B_2$ to $Aut(B_1)$. Since $B_1$ is a cyclic group, $Aut(B_1)$ is a finite group. Hence
up to a subgroup of finite index, we get a central group extension
$$1 \rightarrow B_1 \rightarrow B_2 \rightarrow B_3 \rightarrow 1$$
where $B_1$ is cyclic and $B_3$ is nilpotent. Now it is easy to see $B_2$ is also nilpotent.
\end{proof}
\par Now we use the above lemmas to prove theorem 1.4.
 In fact, assume $\Gamma\simeq \pi_1(M)$, where $M$ is a compact Sasakian manifold. By a theorem of Wadsley [Wa75], we can always assume $M$ is quasi-regular. So we can assume $M$ is an orbifold $S^1$ bundle over some compact K$\ddot{a}$hler orbifold $N$ and $p: M\rightarrow N$ is the projection map. Then we have the following long exact sequence
$$ \rightarrow \pi_2^{orb}(N)\rightarrow \mathbb{Z}\rightarrow \pi_1(M) \rightarrow \pi_1^{orb}(N) \rightarrow 1$$
From which we get the following short exact sequence
$$1 \rightarrow \Gamma_1 \rightarrow \pi_1(M) \rightarrow \pi_1^{orb}(N) \rightarrow 1$$
where $\Gamma_1$ is a cyclic group.
\par Now we prove the first part of theorem 1.4. First note $\Gamma$ can not have $2$ ends otherwise it contains an infinite cyclic subgroup of finite index which is impossible by proposition 1.1. If $\Gamma_1$ is an infinite cyclic group, $\Gamma$ has $1$ end by lemma 3.2. Otherwise $\Gamma_1$ is finite and $\Gamma$ is quasi-isometric to
$\pi_1^{orb}(N)$. It follows that $\Gamma$ has $0$ or $1$ end by lemma 3.1.
\par If $\Gamma$ is solvable, then $\pi_1^{orb}(N)$ is also solvable. It follows that $\pi_1^{orb}(N)$ is a virtually nilpotent group by lemma 3.1 and $\Gamma$ is also virtually nilpotent by lemma 3.3.
This proves the second part of theorem 1.4.
\section*{4. Harmonic maps from compact Sasakian manifolds}
We first state the following theorem due to R. Petit. See [Pet02].
\\
\\
\noindent \textbf{Theorem 4.1}
{\itshape Suppose $f$ is a harmonic map from a compact Sasakian manifold $(M_1, h_1)$ to a compact Riemannian manifold $(M_2, h_2)$ with nonpositive sectional curvature, then $f_*(\xi)=0$, where $\xi$ is the Reeb vector field associated to $M_1$. Moreover, if $M_1$ is a regular Sasakian manifold fibering over K$\ddot{a}$hler manifold $N$ and $p: M_1\rightarrow N$ is the projection map, there exists a harmonic map $g$ from $N$ to $M_2$ such that $f=gp$. }
\\
\par It was known by C.Boyer and K. Galicki [BoG08]
that the connected sum of two compact negatively curved manifolds can not admit Sasakian structure. The following corollary is a generalization of this result.
\\
\\
\noindent \textbf{Corollary 4.2}
{\itshape There is no continous map of nonzero degree from a compact Sasakian manifold to
the connected sum $M_1^n \sharp M_2^n$, where $(M_1^n, g_1)$ is a compact Riemannian manifold with nonpositive sectional curvature and $M_2^n$ is any compact $n$-dimensional manifold.
In particular,  $M_1^n \sharp M_2^n$ can not admit Sasakian structure.}
\begin{proof}
Suppose $f_1$ is a continous map of nonzero degree from a compact Sasakian manifold $(M_0^n, g_0)$ to $M_1^n \sharp M_2^n$. Note
there is a continous map of degree one
$f_2: M_1^n \sharp M_2^n \rightarrow M_1^n$. Let $f_3=f_2f_1$, then $f_3:(M_0^n,g_0)\rightarrow (M_1^n,g_1)$ is a  continous map of nonzero degree.
 By a classic theorem of Eells and Sampson ([ES64]), within the same homotopy class, we can find
a harmonic representive $f_4$. On the other hand, by theorem 4.1, $f_4$ can not be surjective by Sard's theorem.
It follows that the degree of $f_4$ is zero. Contradiction.
\end{proof}
\par
As an application of corollary 4.2, any compact manifold admitting a metric of nonpositive sectional curvature cannot
admit Sasakian structure. For example, suppose $H^3$ is a $3$-dimensional hyperbolic homology sphere and $T^2$ is the $2$-dimensional torus. Then $T^2 \times H^3$ can not admit Sasakian structure by corollary 4.2.
However, as far as the author knows, the nonexistence of Sasakian structure on
$T^2 \times H^3$
does not follow from any previously known obstructions.
\\
\par Now we are about to prove proposition 1.5. For the rest of this section, we always assume $M$ is a compact regular Sasakian manifold fibering over some compact K$\ddot{a}$hler manifold $N$ and $p: M\rightarrow N$ is the projection map. Denote $\Gamma = \pi_1(M)$.
\par
We first show the first part.
Now assume $\Gamma \simeq \pi_1(B)$, where $B=\Gamma \setminus G /K$, where $K$ is a maximal compact subgroup of $G$. It follows that $B$ is a compact locally symmetry space of noncompact type and so admits a Riemannian metric of nonpositive sectional curvature. By a classic theorem of Eells and Sampson, we can choose a harmonic map $f: M\rightarrow B$ inducing the isomorphism between $\pi_1(M)$ and $\pi_1(B)$. By theorem 4.1, there exists a harmonic map $g$ from $N$ to $B$ such that $f=gp$.
It follows that $p_*:\pi_1(M)\rightarrow \pi_1(N)$ is injective.
On the other hand, from the long exact sequence of homotopy groups we know
$p_*$ is surjective and hence $p_*$  is an isomorphism.
Now $\Gamma$ is a K$\ddot{a}$hler group and also a cocompact lattice in $G$, where $G$ is $SO(1,n), n>2$ or $F_{4(-20)}$
or a simple real non-Hermitian Lie group of noncompact type with real rank at least 20. This is impossible by the results of Carlson, Hern$\acute{a}$ndez, Klingler and Toledo. See [CaT89], [CaH91], [Kl10].
\par Now we prove the second part. Suppose $\Gamma =\pi_1(M) \simeq \pi_1(V^3)$, where $V^3$ is a compact $3$-dimensional manifold. By passing to orientable cover, we can assume $V^3$ is orientable. First we have the following long exact sequence of homotopy groups:
$$\rightarrow\pi_2(N) \rightarrow \pi_1(S^1) \rightarrow \pi_1(M) \rightarrow \pi_1(N) \rightarrow 1$$
From it we get a short exact sequence:
$$1 \rightarrow i_*(\pi_1(S^1)) \rightarrow \pi_1(M) \rightarrow \pi_1(N) \rightarrow 1$$
where $i: S^1\rightarrow M$ is the inclusion map. If $i_*(\pi_1(S^1))=\{1\}$, $\Gamma= \pi_1(M) \simeq \pi_1(N)$ and so $\Gamma$ is a K$\ddot{a}$hler group. By assumption, we have $\Gamma \simeq \pi_1(V^3)$.
It follows that $\Gamma$ is finite by
a theorem of A. Dimca and A. Suciu, see [DS09] and [Kot11]. Hence the universal cover of $V^3$ is differemorphic to $S^3$
 by the work of Perelman [Per02], [Per03]. Now we assume $i_*(\pi_1(S^1))$ is a nontrivial cyclic group. Then $\pi_1(V^3)$ contains a nontrivial cyclic normal subgroup. In this case, we first
\\
\\
\noindent \textbf{Claim:} {\itshape $V^3$ is a Seifert manifold.}
\\
\par Given the above claim, it follows that $V^3$ carries one of the geometries $S^2 \times \mathbb{R}, H^2\times \mathbb{R}, {\mathbb{R}}^3, S^3, \widetilde{SL(2,\mathbb{R})}$,
 $Nil$, where $Nil$
is the $3$-dimensional Heisenberg group. But those manifolds carrying one of the geometries
$S^2 \times \mathbb{R}, H^2\times \mathbb{R}, {\mathbb{R}}^3$
 have virtually odd first Betti number. However, any Sasaki group has even first Betti number by proposition 1.1.
 So the universal cover of $V^3$ is differemorphic to $S^3$ or the $3$-dimensional Heisenberg group or $\widetilde{SL(2,\mathbb{R})}$. The other direction is obvious.
 \\
 \\
 \noindent \textbf{Proof of claim }
\par First of all, we can assume $V^3$ is prime in the sense being indecomposable
under connected sum, since a nontrivial free product is never a Sasaki group by theorem
1.4. By assumption, $\pi_1(V^3)$ contains a nontrivial cyclic normal subgroup $\Gamma_0$. Firstly suppose
 $\Gamma_0$ is a finite cyclic group. Then $\pi_1(V^3)$ has nontrivial torsion.
  Since $V^3$ is a prime also orientable 3-manifold, it follows that $\pi_1(V^3)$ must be finite
   by a theorem of Epstein, see [Ep61], [H76]. So the universal cover of $V^3$ is differemorphic to $S^3$
 by the work of Perelman. Now we assume $\Gamma_0$ is an infinite cyclic group.
Since $V^3$ is prime and orientable, it is either $S^2 \times S^1$
or irreducible, that is, every embedded $S^2$ bounds a 3-cell. Because a Sasaki group can not be infinite cyclic by corallary 1.3, we see $V^3$ must be irreducible. On the other hand, $\pi_1(V^3)$ contains an
 infinite cyclic normal subgroup and so
 $V^3$ must be a Seifert manifold by a theorem of A. Casson, D. Jungreis and D. Gabai, see [CaJ94], [Ga92].
\section*{5. A transversal Jacobi equation and proof of theorem 1.6}
The proof of theorem 1.6 is based on the following transversal Cartan-Hadamard theorem:
\\
\\
\noindent \textbf{Theorem 5.1} {\itshape Suppose $(M^{2n+1}, g)$ is a complete Sasakian manifold with nonpositive transversal sectional curvature,
then its universal cover is differemorphic to $\mathbb{R}^{2n+1}.$}
\\
\\
 \noindent \textbf{Corollary 5.2} {\itshape Suppose $(M^{2n+1}, g)$ is a compact Sasakian manifold with nonpositive transversal sectional curvature, then its fundamental group can not be Gromov hyperbolic.}
 \\
 \par In fact, by a theorem of I. Mineyev [Mi01], [Mi02], any compact aspherical manifold has positive simplicial volume
if its fundamental group is Gromov hyperbolic. Since any compact Sasakian manifold has vanishing simplicial volume, corollary 5.2 follows.
\\
\par
By the classification of compact Sasakian manifolds in dimension 3
\\
([BoG08]), we get
\\
\\
\noindent \textbf{Corollary 5.3.} {\itshape Suppose $(M^3, g)$ is a 3-dimensional compact Sasakian manifold with  nonpositive transversal sectional curvature,
then $M^3 $ is differemorphic to the compact quotient of the $3$-dimensional Heisenberg group or $\widetilde{SL(2,\mathbb{R})}$.}
\\
\par On the other hand, compact quotients of the $3$-dimensional Heisenberg group ($\widetilde{SL(2,\mathbb{R})}$) are Seifert circle bundles over flat (hyperbolic) orbifolds and so admit Sasakian structure with nonpositive  transversal  sectional curvature.
\par Before we prove theorem 5.1, we show how to derive theorem 1.6 from it. We prove it by contradiction.
 Suppose $\Gamma=\pi_1(M,g)\simeq\pi_1(N,h)$, where $(M,g)$ is a compact Sasakian manifold with nonpositive transversal sectional curvature and $(N,h)$ is a compact manifold with nonpositive sectional curvature.
 Then by a classic result of Eells and Sampson, there exists a harmonic map $f$ which induces the isomorphism of the fundamental groups. By a theorem of A. Banyaga and P. Rukimbira ([Ban90], [Ruk95], [Ruk99]), there exists a closed leaf of the characteristic foliation on $M$. Since $M$ is compact, we see this leaf is differemorphic to $S^1$. From the proof of theorem 5.1, we see this $S^1$ is lifted to a real line under covering map
$p:\widetilde{M}\rightarrow M$, where $\widetilde{M}$ is the universal cover of $M$.
 This implies the inclusion
map $i_*: \pi_1(S^1)\rightarrow \pi_1(M)$ is injective. On the other hand, since $f$ is harmonic, using theorem 4.1,
$f_*(\xi)=0$, where $\xi$ is the Reeb vector field associated to $M$. So $f$ maps
a nontrivial subgroup of $\Gamma$ to zero, which contradicts that $f$ is an isomorphism.
\\
\par Now we are in position to prove theorem 5.1.
 The idea is similar to the proof of Cartan-Hadamard theorem. Suppose $\widetilde{M^{2n+1}}$  is the universal cover of $M^{2n+1}$ .
 Then it is also a complete Sasakian manifold with nonpositive transversal sectional curvature.
Suppose $\alpha(s)$ is a leaf of the characteristic foliation $F_{\xi}$ on $\widetilde{M^{2n+1}}$.
We show its normal exponential map $exp^\bot: \alpha(s)^\bot\rightarrow \widetilde{M^{2n+1}} $ is a differmorphism, where
$\alpha(s)^\bot$
is the normal bundle of $\alpha(s)$.
To do this, suppose $\gamma(t)$, $t\in[0,1]$ is any minimizer geodesic which is perpendicular to
 $\alpha(s)$ at $\alpha(0)=\gamma(0)$
and $J(t)$ is any Jacobi field
along $\gamma(t)$ such that $J(0)=\lambda \xi(0)$ and $J(1)=0$, where $\xi(0)=\alpha'(0)$. It suffices to show $J(t)\equiv 0$.
Then we see there is no focal point to $\alpha(s)$ and so $exp^\bot$ is a differmorphism.
\par There is a natural splitting $T_{\gamma(t)}M= \xi(t)\oplus \xi(t)^\bot$, where $\xi(t)$ is the tangential part
to leaves of the characteristic foliation $F_{\xi}$ and $\xi(t)^\bot$ is the orthogonal part.
Following B.Wilking [Wi07], define $Y(t)=J^\bot(t)$ and
 $\nabla^\bot_{\frac{\partial}{\partial t}}Y=(\nabla_{\frac{\partial}{\partial t}}Y)^\bot=(\nabla_{\dot{\gamma}}Y)^\bot$, where we denote
 $X^\bot(t)$ is the projection of $X(t)$ to $\xi(t)^\bot$ for any vector field $X(t)$ along $\gamma(t)$.
 \\ \noindent \textbf{Lemma 5.4} $Y(t)$
 satisfies the following transversal Jacobi equation:
 $$\nabla^\bot_{\frac{\partial}{\partial t}}\nabla^\bot_{\frac{\partial}{\partial t}}Y+ R^T(Y, \dot{\gamma}) \dot{\gamma}=0$$
 Before we prove it, we need the following lemma.
  \\ \noindent \textbf{Lemma 5.5}
 $$ R^T(X,Y)Y=(R(X,Y)Y)^\bot+3<X,\phi(Y)>\phi(Y), where X,Y\perp \xi.$$
 \begin{proof}
 Choose any vector field $Z$ such that $Z\bot\xi$. We see
 $$ R^T(X,Y,Y,Z)=<R^T(X,Y)Y,Z>=<\nabla_X^T\nabla_Y^TY-\nabla_Y^T\nabla_X^TY-\nabla^T_{[X,Y]}Y,Z>$$
$$ =<\nabla_X(\nabla_YY-<\nabla_YY,\xi>\xi)-\nabla_Y(\nabla_XY-<\nabla_XY,\xi>\xi)$$
$$-\nabla^T_{[X,Y]-<[X,Y],\xi>\xi}Y-\nabla^T_{<[X,Y],\xi>\xi}Y,Z>$$
$$=<R(X,Y)Y,Z>+<\nabla_XY,\xi><\nabla_Y\xi,Z>+<[X,Y],\xi><\nabla_Y\xi,Z>$$
$$=<R(X,Y)Y,Z>+3<X,\phi(Y)><\phi(Y),Z>$$
where the last equality follows from $X, Y \perp \xi$ and $\xi$ is a Killing vector field.
\end{proof}
Now we are in position to prove lemma 5.4. It suffices to prove it at generic $t_0$, i.e. $J(t_0)\neq0$.
First note we can assume $Y(t_0)=J(t_0)$.
Choose vector fileds $X_i(t), i=1,2,\cdots,2n$ such that
$X_i\bot\xi$ and also $J(t_0)=X_1(t_0), \nabla^\bot_{\frac{\partial}{\partial t}}X_i(t)=0, <X_i,X_j>(t)=\delta_{ij}$ for all $t\in[0,1]$. First note we have
$$X_i'=<X_i',\xi>\xi=<\nabla_{\dot{\gamma}}X_i,\xi>\xi=-<X_i,\nabla_{\dot{\gamma}}\xi>\xi$$
Using this, at $t_0$, we have
$$<\nabla^\bot_{\frac{\partial}{\partial t}}\nabla^\bot_{\frac{\partial}{\partial t}}Y, X_i>
=<\nabla_{\frac{\partial}{\partial t}}\nabla^\bot_{\frac{\partial}{\partial t}}Y, X_i>$$
$$
=\frac{\partial}{\partial t}<\nabla^\bot_{\frac{\partial}{\partial t}}Y,X_i>-<\nabla^\bot_{\frac{\partial}{\partial t}}Y,\nabla^\bot_{\frac{\partial}{\partial t}}X_i>$$
$$=\frac{\partial}{\partial t}<\nabla^\bot_{\frac{\partial}{\partial t}}Y,X_i>
=\frac{\partial^2}{\partial t^2}<J,X_i>$$
$$
=<J'',X_i>+2<J',X_i'>+<J,X_i''>$$
$$=-<R(Y, \dot{\gamma}) \dot{\gamma}, X_i>+2<J',X_i'>+<X_1,X_i''>$$
$$=-<R(Y, \dot{\gamma}) \dot{\gamma}, X_i>+2<J',X_i'>$$
$$+
\frac{d}{dt}_{t=t_0}<X_1,X_i'>-<X_1',X_i'>$$
$$=-<R(Y, \dot{\gamma}) \dot{\gamma}, X_i>+2<J',X_i'>-<X_1',X_i'>$$
$$=-<R(Y, \dot{\gamma}) \dot{\gamma}, X_i>-2<J',\xi><X_i,\xi'>-<X_1,\xi'> <X_i,\xi'>$$
$$=-<R(Y, \dot{\gamma}) \dot{\gamma}, X_i>-2<J,\xi'><X_i,\xi'>
-<J,\xi'> <X_i,\xi'>$$
$$=<-R^T(Y, \dot{\gamma}) \dot{\gamma},X_i>$$
 where the last equality follows from lemma 5.5
and the last second equality follows from
 \\
 \noindent \textbf{Lemma 5.6}
 $$<J'(t), \xi(t)>=<J(t),\xi'(t)>$$
 \begin{proof}
$$<J'(t), \xi(t)>=<\nabla_{\dot{\gamma}}J,\xi>(t)$$
$$=<\nabla_J\dot{\gamma},\xi>(t)
=-<\dot{\gamma},\nabla_J\xi>(t)$$
$$=<J,\nabla_{\dot{\gamma}}\xi>(t)=<J(t),\xi'(t)>$$
 \end{proof}
 \par Now define $f(t)=\frac{1}{2}\|Y(t)\|^2$, then $f(0)=f(1)=0$ since $Y(0)=Y(1)=0$. Moreover, we have
 $$f'(t)=<\nabla^\bot_{\frac{\partial}{\partial t}}Y,Y>$$
 and also
 $$f''(t)=<\nabla^\bot_{\frac{\partial}{\partial t}}Y,\nabla^\bot_{\frac{\partial}{\partial t}}Y>+<\nabla^\bot_{\frac{\partial}{\partial t}}\nabla^\bot_{\frac{\partial}{\partial t}}Y,Y>$$
$$ = <\nabla^\bot_{\frac{\partial}{\partial t}}Y,\nabla^\bot_{\frac{\partial}{\partial t}}Y>-R^T(Y, \dot{\gamma}, \dot{\gamma},Y)\geq 0$$
From this we know $f(t)$ is a convex function with $f(0)=f(1)=0$ and so $f(t)=0$ for any $t\in[0,1]$.
So $Y(t)\equiv0$ and we can assume $J(t)=\lambda(t)\xi(t)$ with $\lambda(0)=\lambda$. Then $\lambda(t)=<J(t),\xi(t)>$ and $\lambda(1)=0$
since $J(1)=0$.
 Taking derivative with respect to $t$, we get
  $$\lambda'(t)=<J'(t), \xi(t)>+<J(t),\xi'(t)>$$
 $$=2<J(t),\xi'(t)>=2<\lambda(t)\xi(t),\xi'(t)>=0.$$
so $J(t)\equiv0$ and $exp^\bot: \alpha(s)^\bot\rightarrow \widetilde{M^{2n+1}}$ is a differmorphism.
Hence $\widetilde{M^{2n+1}}$ is diferemorphic to the normal bundle of a 1-dimensional manifold.
Since $\widetilde{M^{2n+1}}$ is simply connected, we see $\widetilde{M^{2n+1}}$ is differemorphic to $\mathbb{R}^{2n+1}.$
\section*{6. Abel-Jacobi maps of compact Sasakian manifolds}
 We first recall the construction of Abel-Jacobi map of any compact Riemannian manifold, see [Ka07] for more details.
 \par Suppose $(M,g)$ is a compact Riemannian manifold. Let $\pi=\pi_1(M)$ be its fundamental group.
 Set $f: \pi \rightarrow \pi^{ab} $ be the abelianisation map and $g: \pi^{ab}\rightarrow \pi^{ab}/tor$
 be the quotient by torsion. Suppose $\bar{M}$ is the covering space of $M$ with $\pi_1(\bar{M})=ker (\phi)$, where $\phi=gf$.
\par Let $E$ be the space of harmonic $1$-form on $M$, with dual $E^*$ canonically identified with $H_1(M, \mathbb{R})$. Fix a basepoint $x_0 \in M$. Then any point $x$ in the universal cover $\tilde{M}$ of $M$ is represented by a point of $M$ together with a path $c$ from $x_0$ to it. By integrating along the path $c$, we get a linear form, $h \rightarrow \int_c h$, on $E$. We thus obtain a map $\tilde{M}\rightarrow E^*=H_1(M, \mathbb{R})$, which descends to a map
 $$\bar{A_M}: \bar{M}\rightarrow E^*, c\mapsto (h\mapsto  \int_c h),$$
 \par By definition, the Jacobi torus of $M$ is the torus
 $$J_1(M)=H_1(M, \mathbb{R})/{H_1(M, \mathbb{Z})_\mathbb{R}}$$
 and the Abel-Jacobi map
 $$A_M: M\rightarrow J_1(M),$$
 is obtained from the map $\bar{A_M}$ by passing to quotients. From the construction,
 it is not hard to see the Abel-Jacobi map induces an isomorphism between the first homology groups with real coefficient.
 \\
 \\
 \noindent \textbf{Proposition 6.1} {\itshape Suppose $(M,g)$ is a compact Sasakian manifold and
 $F_{\xi}$ is the characteristic foliation on $M$. Then the restriction of the Abel-Jacobi map of $M$ to
 any leaf of $F_{\xi}$ is a constant map.}
 \begin{proof}
 By a theorem of Tachibana ([Tac65]), any harmonic one form $h$ on $M$ satisfies $h(\xi)=0$. Using it, proposition 6.1 easily follows from the construction of Abel-Jacobi map.
 \end{proof}
 \noindent \textbf{Corollary 6.2} {\itshape Suppose $(M,g)$ is a compact Sasakian manifold and $i: S^1\rightarrow M$ is the inclusion map, where $S^1$ is a closed leaf of the characteristic foliation
 $F_{\xi}$ on $M$. Then
 $i_*(\pi_1(S^1))\subseteq ker \phi$, where $\phi$ is the map constructed
 in the beginning of this section. In other words, $S^1$ generates a trivial element in $H_1(M, \mathbb{R})$.}
 \begin{proof}
  By proposition 6.1, the Abel-Jacobi map of $M$ takes $S^1$ to a point. However, the Abel-Jacobi map induces an isomorphism between the first homology groups with real coefficient and so $S^1$ generates a trivial element in $H_1(M, \mathbb{R})$.
\end{proof}

\bigskip

\noindent
{
Xiaoyang Chen \\
{\small Department of Mathematics}\\
{\small University of Notre Dame} \\
{\small 46556,Notre Dame, Indiana, USA}\\
{\small 574-631-3741} \\
{\small xchen3@nd.edu}

\begin{thebibliography}{999999}
\bibitem[ABCKT96]{ABCKT96} J. Amor$\acute{o}$s, M. Burger, K. Corlette, D. Kotschick and D. Toledo, Fundamental groups
of compact K$\ddot{a}$hler manifolds, Mathematical Suveys and Monographs, vol. 44, Amer. Math. Soc., Providence, R.I. 1996.
\bibitem[Ban90]{Ban90} A. Banyaga, A note on Weinstein's conjecture, Proc. Amer. Math. Soc. 109 (1990), no 3, 855-858.
\bibitem[Be68]{Be68} G. Bergman, On groups acting on locally finite graphs, Ann. of Math. 88 (1968), 335-340.
\bibitem[Bl02]{Bl02}D. E. Blair, Riemannian geometry of contact and symplectic manifolds, Progress in mathematics, 203, Birkh\"{a}user, Boston, 2002.
\bibitem[BG67]{BG67}D. E. Blair and S. I. Goldberg, Topology of almost contact manifolds,  J. Differential Geom.
1 (1967), 347-354.
\bibitem[BoG99]{BoG99} C. Boyer and K. Galicki, 3-Sasakian manifolds, in: C.LeBrun and M. Wang. (eds), Surveys in Differetial Geometry, Vol. VI: Essays on
Einstein manifolds, International Press, Boston, MA, 1999, 123-184.
\bibitem[BoG08]{BoG08} C. Boyer and K. Galicki, Sasakian geometry, Oxford University Press (2008).
\bibitem[C11]{C11} F. Campana, Quotients résoluble ou nilpotents des groupes
de K$\ddot{a}$hler orbifoldes,  Manuscripta Math. 135 (2011), no. 1-2, 117-150.
\bibitem[CaH91]{CaH91} J. A. Carlson and L. Hern$\acute{a}$ndez, Harmonic maps from compact K$\ddot{a}$hler manifolds to exceptional hyperbolic spaces,  J. Geom. Anal. 1 (1991), 339-357.
\bibitem[CaJ94]{CaJ94} A. Casson and D. Jungreis, Convergence groups and Seifert fibered 3-manifolds, Invent. Math. 118 (1994), 441-456.
\bibitem[CaT89]{CaT89} J. A. Carlson and D. Toledo, Harmonic mappings of K$\ddot{a}$hler manifolds to locally symmetry spaces, Publ. Math. I.H.E.S, 69 (1989), 173-201.
\bibitem[CaT95]{CaT95} J. A. Carlson and D. Toledo, Quadratic presentations and nilpotent K$\ddot{a}$hler
groups, J. Geom. Anal. 5 (1995), 359-377.
\bibitem[Co70]{Co70} D. E. Cohen, Ends and free products of groups, Math. Z. 114 (1970), 9-18.
\bibitem[Da08]{Da08} M. W. Davis, Lectures on orbifolds and reflection groups, Lectures in the summer school on transformations groups and orbifolds, Zhejiang University, Hangzhou, China, 2008.
\bibitem[Dea04]{Dea04}O. Dearricott, Positive sectional curvature on $3$-Sasakian manifolds, Ann. Global Anal. Geom. 25(2004), no. 1, 59-72.
\bibitem[Del10]{Del10} T. Delzant, L'invariant de Bieri Neumann Strebel des groupes fondamentaux des vari$\acute{e}$t$\acute{e}$s k$\ddot{a}$hl$\acute{e}$riennes, Math. Annalen 348 (2010), 119-125.
\bibitem[DK]{DK} C. Drutu and M. Kapovich, Lectures on geometric group theory. http://www.math.ucdavis.edu/~kapovich/EPR/ggt.pdf.
\bibitem[DS09]{DS09} A. Dimca and A. Suciu, Which $3$-manifold groups are K$\ddot{a}$hler groups?  J. Eur. Math. Soc. 11 (2009), 521-528.
\bibitem[Eb82]{Eb82} P. Eberlein, A Canonical form for compact nonpositively curved manifolds
whose fundamental groups have nontrivial center, Math. Ann. 260 (1982), 23-29.
\bibitem[Ep61]{Ep61} D. B. A. Epstein, Projective planes in 3-manifolds, Proc. London Math. Soc. (3) 11 (1961) 469-484.
\bibitem[ES64]{ES64} J. Eells and J. H. Sampson, Harmonic maps of Riemannian manifolds, Amer. Jour. Math.
86 (1964), 109-160.
\bibitem[F66]{F66} T. Fujitani, Complex-valued differential forms on normal contact Riemannian manifolds,
T\^{o}hoku Math. J. 18 (1966), 349-361.
\bibitem[FOW09]{FOW09} A. Futaki, H. Ono and G. Wang,
 Transverse K$\ddot{a}$hler geometry of Sasaki manifolds and toric Sasaki-Einstein manifolds,
 \\  J. Differential Geom. 83 (2009), no. 3, 585-635.
\bibitem[Ga92]{Ga92} D. Gabai, Convergence groups are Fuchsian group, Ann. of Math. 136 (1992),447-510.
\bibitem[Gr89]{Gr89} M. Gromov, Sur le groupe fondamental d'une vari$\acute{e}$t$\acute{e}$ k$\ddot{a}$hl$\acute{e}$rienne,
C. R. Acad. Sci. Paris S$\acute{e}$r. I Math. 308 (1989), 67-70.
\bibitem[H76]{H76} J. Hempel, 3-manifolds,  Ann. Math. Stud. no 86, Princeton University Press (1976).
\bibitem[HW94]{HW94} W. Heil and W. Whitten, The Seifert fiber space conjecture and torus theorem for nonorientable $3$-manifolds, Canad. Math. Bull. 37 (1994), 482-489.
\bibitem[Ka07]{Ka07} M. G. Katz, Systolic geometry and topology, Mathematical Surveys and Monographs, vol 137, American Mathematical Society, 2007.
\bibitem[Kl10]{Kl10} B. Klingler, On the cohomology of K$\ddot{a}$hler groups, arxiv:1005. 2835v1 [math. GR] 17 May 2010.
\bibitem[Kob63]{Kob63} S. Kobayashi, Topology of positively pinched K$\ddot{a}$hler manifolds, T\^{o}hoku Math. J. 15 (1963), 121-139.
\bibitem[Kot11]{Kot11} D. Kotschick, Three manifolds and K$\ddot{a}$hler groups. arXiv:1011.4084v3 [math. GT] 12 Feb 2011.
\bibitem[Le95]{Le95} B. Leeb, $3$-manifolds with(out) metrics of nonpositive curvature, Invent. Math. 122 (1995) 277-289.
\bibitem[Lu88]{Lu88} L. Luecke, Finite covers of $3$-manifolds containing essential tori, Trans. Amer. Math. Soc. 310 (1988), 381-391.
\bibitem[Mi01]{Mi01} I. Mineyev, Straightening and bounded cohomology of hyperbolic groups, Geom. Funct. Anal. 11 (2001), no 4, 807-839.
\bibitem[Mi02]{Mi02} I. Mineyev, Bounded cohomology characterize hyperbolic groups, Q. J. Math. 53 (2002),
no 1, 59-73.
\bibitem[MSY08]{MSY08} D. Martelli, J. Sparks and S. T. Yau, Sasaki-Einstein manifolds and volume minimisation, Commun. Math. Phys. 280 (2008), 611-673.
\bibitem[Per02]{Per02} G. Perelman, The entropy formula for the Ricci flow and its geometric applications, Preprint, arXiv: math/0211159v1 [math.DG] 11 Nov 2002.
\bibitem[Per03]{Per03} G. Perelman, Ricci flow with surgery on three manifolds, Preprint, arXiv: math/0303109v1 [math.DG] 10 Mar 2003.
\bibitem[Pet02]{Pet02} R. Petit, Harmonic maps and strictly pseudoconvex CR manifolds, Comm. Anal. Geom. 10 (2002), no 3, 575-610.
\bibitem[Ruk95]{Ruk95} P. Rukimbira, Topology and closed characteristics of K-contact manifolds, Bull. Belg. Math. Soc. Simon Stevin, 2 (1995), no. 3, 349-356.
\bibitem[Ruk99]{Ruk99} P. Rukimbira, On K-contact manifolds with minimal number of closed characteristics,
Proc. Amer. Math. Soc. 127 (1999) no. 11, 3345-3351.
\bibitem[Sc83]{Sc83} P. Scott, The geometries of $3$-manifolds, Bull. London Math Soc. 15 (1983), 401-487.
\bibitem[St68]{St68} J. R. Stalling, On torsion-free groups with infinitely many ends, Ann. of Math. 88 (1968), 312-334.
\bibitem[St71]{St71} J. R. Stalling, Group theory and three dimensional manifolds, vol 4, Yale Mathematical Monographs, Yale University Press, New Haven, Conn-London, 1971.
\bibitem[Tac65]{Tac65} S. Tachibana, On harmonic tensors in compact Sasakian spaces,  T\^{o}hoku Math. J.,
17 (1965), 271-284.
\bibitem[Tan67]{Tan67} S. Tanno, Harmonic forms and Betti numbers of certain contact Rimannian manifolds, J. Math. Soc. Japan. 19 (1967) 308-316.
\bibitem[Th78]{Th78} W. Thurston, The geometry and topology of three-manifolds, Princeton University lectures notes, 1978-1981.
\bibitem [Wa75]{Wa75} A. W. Wadsley, Geodesic foliations by circles, J. Differential Geom. 10 (1975), no. 4, 541-549.
\bibitem [Wi07]{Wi07} B. Wilking, A duality theorem for Riemannian foliations in nonnegative sectional curvature,
GAFA. 17 (2007), 1297-1320.
\end{thebibliography}
\end{document}